\documentclass[12pt]{amsart}
\usepackage{amssymb}
\usepackage{verbatim}
\usepackage{amscd}

\theoremstyle{plain}
\newtheorem{theorem}{Theorem}
\newtheorem{corollary}{Corollary}
\newtheorem{proposition}{Proposition}
\newtheorem{lemma}{Lemma}

\theoremstyle{definition}
\newtheorem{definition}{Definition}

\theoremstyle{remark}
\newtheorem{remark}{Remark}

\newtheorem*{notation}{Notation}

\newcommand{\ti}{\tilde}

\def\bbC{\mathbb C}

\def\bbN{\mathbb N}

\def\bbT{\mathbb T}
\def\bbZ{\mathbb Z}

     \newcommand{\sA}{\mathcal A}
     \newcommand{\sB}{\mathcal B}
     \newcommand{\sC}{\mathcal C}
     
     \newcommand{\sE}{\mathcal E}

     \newcommand{\sH}{\mathcal H}

\newcommand{\fO}{\mathfrak O}     \newcommand{\sO}{\mathcal O}

     \newcommand{\sS}{\mathcal S}
\newcommand{\fT}{\mathfrak T}

     \newcommand{\sY}{\mathcal Y}

\newcommand{\al}{\alpha}
\newcommand{\be}{\beta}
\newcommand{\vpi}{\varphi}
\newcommand{\si}{\sigma}
\newcommand{\la}{\lambda}

\begin{document}

\title[C$^*$-envelope and Nest Representations]
{The C$^*$-envelope of a semicrossed product and Nest Representations}

\author{Justin R. Peters}
\address{Department of Mathematics\\
    Iowa State University, Ames, Iowa, USA and}

\email{peters@iastate.edu}


 \begin{abstract} Let $X$ be compact Hausdorff, and $\vpi: X \
 \to X$ a continuous surjection. Let $\sA$ be the semicrossed product
 algebra corresponding to the relation $fU = Uf\circ \vpi$. Then the
 C$^*$-envelope of $\sA$ is the crossed product of a commutative
 C$^*$-algebra which contains $C(X)$ as a subalgebra, with respect
 to a homeomorphism which we construct. We also show there
 are``sufficiently many'' nest representations.
 \end{abstract}
\maketitle

\section{Introduction} \label{s:intro} In \cite{Pet.Sem} the notion of the
semi-crossed product of a C$^*$-algebra with respect to an
endomorphism was introduced. This agreed with the notion of a
nonselfadjoint or analytic crossed product introduced earlier by
McAsey and Muhly (\cite{MM.Rep}) in the case the endomorphism was an
automorphism. Neither of those early papers dealt with the
fundamental question of describing the C$^*$-envelopes of the class
of operator algebras being considered.

That open question was breached in the paper \cite{MS.Ten}, in which
Muhly and Solel described the C$^*$-envelope of a semicrossed
product in terms of C$^*$-correspondences, and indeed determined the
C$^*$-envelopes of many classes of nonselfadjoint operator algebras.

While it is not our intention to revisit the results of
\cite{MS.Ten} in any detail, we recall briefly what was done.  Given
a C$^*$-algebra $\sC$ and an endomorphism $\al$ of $\sC$ one forms
the semicrossed product $ \sA := \sC\rtimes_{\al} \bbZ^+ $ as
described in Section~\ref{s:semi}.  First one views $\sC$ as a
C$^*$-correspondence $\sE$ by taking $\sE = \sC$ as a right $\sC$
module, and the left action given by the endomorphism. One then
identifies the tensor algebra (also called the analytic Toeplitz
algebra) $\fT_+(\sE)$ with the semicrossed product $\sA$. The
C$^*$-envelope of $\sA$ is given by the Cuntz-Pimsner algebra
$\fO(\sE)$.

The question that motivated this paper was to find the relation
between the C$^*$-envelopes of semicrossed products, and crossed
products. Specifically, when is the C$^*$-envelope of a semicrossed
product a crossed product? If the endomorphism $\al$ of $\sC$ is
actually an automorphism, then the crossed product
$\sC\rtimes_{\al}\bbZ$ is a natural candidate for the
C$^*$-envelope, and indeed, as noted in \cite{MS.Ten}, this is the
case. In this paper we answer that question in case the
C$^*$-algebra $\sC$ is commutative (and unital).  Indeed, it turns
out that the C$^*$-envelope is \textit{always} a crossed product (cf
Theorem~\ref{t:maint}).

For certain classes of nonselfadjoint operator algebras, nest
representations play a fundamental role akin to that of the
irreducible representations in the theory of C$^*$-algebras.  The
notion of nest representation was introduced by Lamoureux
(\cite{Lam.Nest}, \cite{Lam.Ide}) in a context with similarities to
that here.  We do not answer the basic question as to whether nest
representations suffice for the kernel-hull topology; i.e., every
closed ideal in a semicrossed product is the intersection of the
kernels of the nest representations containing it. What we do show
is that nest representations suffice for the norm: the norm of an
element is the supremum of the norms of the isometric covariant nest
representations (Theorem~\ref{t:nestrep}). The results on nest
representation require some results in topological dynamics, which,
though not deep, appear to be new.

The history of work in anaylytic crossed products and semicrosed
products goes back nearly forty years. While in this note we do not
review the literature of the subject, we mention the important paper
\cite{DKM.Jac} in which the Jacobson radical of a semicrossed
product is determined and necessary and sufficient conditions for
semi-simplicity of the crossed product are obtained. We use this in
Proposition~\ref{p:semis} to show that the simplicity of the
C$^*$-envelope implies the semisimplicity of the semicrossed
product.

 In very recent work of Davidson and Katsoulis
(\cite{DK.Con}), semicrossed products are viewed as an example of a
more general class of Banach Algebras associated with dynamical
systems which they call conjugacy algebras. They have extracted
fundamental properties needed to obtain, for instance, the result
that conjugacy of dynamical systems is equivalent to isomorphism of
the conjugacy algebras. It would be worthwhile to extend the results
here to the broader context.

\section{Dynamical Systems} \label{s:prel}  In our context, $X$ will denote a compact
Hausdorff space.  By a \textit{dynamical system} we will simply mean
a space $X$ together with a mapping $\vpi : X \to X$. In this
article, the map $\vpi$ will always be a continuous surjection.

\begin{definition} \label{d:ext}
Given a dynamical system $(X, \vpi)$ we will say (following the
terminology of \cite{PP.Dyn}) the dynamical system $(Y, \psi)$ is an
\emph{extension} of $(X, \vpi)$ in case there is a continuous
surjection $p: Y \to X$ such that the the diagram
\[
\begin{CD}
Y        @>\psi>>       Y \\
@VpVV       @VpVV \\
X         @>\vpi>>       X \\
\end{CD} \qquad (\dag)
\]
commutes. The map $p$ is called the extension map (of $Y$ over $X$).
\end{definition}

\begin{notation} In case $p$ is a homeomorphism, it is called a
\emph{conjugacy}.
\end{notation}

Given a dynamical system $(X, \vpi)$ there is a canonical procedure
for producing an extension $(Y, \psi)$ in which $\psi$ is a
homeomorphism.

Let $\ti{X} = \{ (x_1, x_2, \dots ):\ x_n \in X \text{ and } x_{n} =
\vpi(x_{n+1}),\ n = 1, 2, \dots \} .$  As $\ti{X}$ is a closed
subset of the product $\Pi_{n=1}^{\infty} X_n$ where $X_n = X,\ n=1,
2, \dots$, so $ \ti{X}$ is compact Hausdorff. Define a map
$\ti{\vpi}: \ti{X} \to \ti{X}$ by
\[ \ti{\vpi}(x_1, x_2, \dots) = (\vpi(x_1), x_1, x_2, \dots) .\]
This is continuous, and has an inverse given by
\[ \ti{\vpi}^{-1}(x_1, x_2, \dots) = (x_2, x_3, \dots) .\]

Define a continuous surjection $p: \ti{X} \to X$ by
\[ p(x_1, x_2, \dots) = x_1. \]
With the map $p$, the system $(\ti{X}, \ti{\vpi})$ is an extension
of the dynamical system $(X, \vpi)$ in which the dynamics of the
extension is given by a homeomorphism.

\begin{definition} \label{d:homeoext} In the case of an extension
 in which the dynamics is given by
a homeomorphism, we will say the extension is a \emph{homeomorphism
extension}.
\end{definition}

\begin{notation} We will call the extension $(\ti{X}, \ti{\vpi})$ the
\emph{canonical} homeomorphism extension.  If $\ti{x} \in \ti{X},\
\ti{x} = (x_1, x_2, \dots)$, we will say that $(x_1, x_2, \dots)$
are the \emph{coordinates} of $\ti{x}$.
\end{notation}

\begin{definition} Given a dynamical system $(X, \vpi)$, a
homeomorphism extension $(Y, \psi)$ is said to be \emph{minimal} if,
whenever $(Z, \si)$ has the property that it is a homeomorphism
extension of $(X, \vpi)$, and $(Y, \psi)$ is an extension of $(Z,
\si)$ such that the composition of the extension maps of $Z$ over
$X$ with the extension map of $Y$ over $Z$ is the extension map of
$Y$ over $X$, then $(Y, \psi)$ and $(Z, \si)$ are conjugate.
\end{definition}

\begin{lemma} \label{l:minext}
Let $(X, \vpi)$ be a dynamical system.  Then the
canonical homeomorphism extension $(\ti{X}, \ti{\vpi})$ is  minimal.
\end{lemma}

\begin{proof} Suppose $(Z, \si)$ is a homeomorphism extension of
$(X, \vpi), \ p: \ti{X} \to Z$ and $q: Z \to X$ are continuous
surjections, and the diagram
\[
\begin{CD}
\ti{X}  @>\ti{\vpi}>>    \ti{X} \\
    @VpVV       @VpVV   \\
Z       @>\si>>     Z   \\
    @VqVV       @VqVV   \\
X       @>\vpi>>    X
\end{CD}
\]
commutes and the composition $q \circ p$ is the extension map of
$\ti{X}$ over $X$, i.e., the projection onto the first coordinate.

Observe that the canonical homeomorphism extension $(\ti{Z},
\ti{\psi})$ of $(Z, \psi)$ is in fact conjugate to $(Z, \psi)$.
Indeed, the map $ z \in Z \mapsto (z, \psi^{-1}(z), \psi^{-2}(z),
\dots )$ is a conjugacy.  Thus it is enough to show that $(\ti{X},
\ti{\vpi})$ is conjugate to $(\ti{Z}, \ti{\si}).$

 Define
a map $r: \ti{Z} \to \ti{X}$ by
\[ \ti{z} := (z, \si^{-1}(z), \si^{-2}(z), \dots) \in \ti{Z} \mapsto
\ti{x}:= (q(z), q(\si^{-1}(z), q(\si^{-2}(z)), \dots) .\]
Observe that this maps into $\ti{X}$, since
\begin{align*}
\vpi(q(\si^{-(n+1)}(z))) &= q(\si(\si^{-(n+1)}(z)))   \\
&= q(\si^{-n}(z))
\end{align*}

Next we claim $r$ maps onto $\ti{X}$.  Let $\ti{x} = (x_1, x_2,
\dots )$ be any element of $\ti{X}$.  Let $z_n \in Z$ be any element
such that $q(z_n) = x_n, \ n= 1, 2, \dots.$ Let
\[ \ti{z_n} := (\si^{n-1}(z_n),\dots, \underset{n}{z_n},
\si^{-1}(z_n), \dots ) .\]

 A subsequence of $\{ \ti{z_n} \}$ converges, say,
to $\ti{z}$.  Since $r(\ti{z}_m)$ agrees with $\ti{x}$ in the first
$n$ coordinates for all $m \geq n$, it follows that $r(\ti{z}) =
\ti{x}.$

To show that $r$ is one-to-one, define a map $\ti{p}: \ti{X} \to
\ti{Z}$ by
\[ \ti{p}(\ti{x}) = (p(\ti{x}), \si^{-1}\circ p(\ti{x}),
\si^{-2}\circ p(\ti{x}) \dots ).\] Note that the fact that $p:
\ti{X} \to Z$ is surjective implies that $\ti{p}$ is surjective. Let
$\ti{x} = (x_1, x_2, x_3, \dots ) \in \ti{X}$.  Then

\begin{align*}
r \circ \ti{p}(\ti{x}) &= r(p(\ti{x}),\ \si^{-1}\circ p (\ti{x}),\
            \si^{-2}\circ p(\ti{x}), \dots ) \\
            &= (q\circ p(\ti{x}),\ q\circ \si^{-1} \circ p(\ti{x}),\ q
            \circ \si^{-2} \circ p(\ti{x}), \dots ) \\
            &= (x_1,\ q\circ p \circ \ti{\vpi}^{-1}(\ti{x}),\ q \circ
            p \circ \ti{\vpi}^{-2}(\ti{x}), \dots ) \\
            &= (x_1, x_2, x_3, \dots ) = \ti{x}.
\end{align*}

where we have used the fact that $q \circ p$ is the projection onto
the first coordinate of $\ti{x}$.  Since $\ti{p}$ is surjective and
$r \circ \ti{p}$ is injective, it follows that $r$ is injective, and
hence $r$ is a conjugacy.

\end{proof}

\begin{lemma}  \label{l:unique}
Let $(X, \vpi)$ be a dynamical system, and let
$(Y, \psi)$ be a minimal homeomorphism extension.  Then $(Y, \psi)$
is conjugate to the canonical homeomorphism extension, $(\ti{X},
\ti{\vpi})$.
\end{lemma}

\begin{proof} By assumption there is a continuous surjection $p: Y
\to X$ such that the diagram
\[
\begin{CD}
Y @>\psi>>  Y  \\
@VqVV       @VqVV \\
X   @>\vpi>>    X
\end{CD}
\]
commutes.

Consider the diagram
\[
\begin{CD}
Y   @>\psi>>    Y \\
@V\ti{q}VV      @V\ti{q}VV \\
\ti{X} @>\ti{\vpi}>>   \ti{X} \\
@VpVV           @VpVV \\
X       @>\vpi>>    X
\end{CD}
\]
where $p$ denotes the canonical extension map of $\ti{X}$ over $X$
(i.e., the projection onto the first coordinate), and the map
$\ti{q}$ is defined as follows:
\[ \text{For } y \in Y,\ \ti{q}(y) = (q(y), q\circ \psi^{-1}(y),
q\circ \psi^{-2}(y), \dots). \] Note that the image lies in $\ti{X}$
since $\vpi(q\circ \psi^{-(n+1)}(y) = q\circ \psi \circ
\psi^{-(n+1)}(y) =q\circ \psi^{-n}(y)$.

Next, observe that $p \circ \ti{q}(y) = q(y)$, so the extension
property is satisfied.  Hence, by definition of minimality of the
homeomorphism extension $(Y, \psi)$, the map $\ti{q}$ is a
conjugacy.

\end{proof}

\begin{corollary} \label{c:uniqminext}
Let $(X, \vpi)$ be a dynamical system.  Then there exists a minimal
homeomorphism extension $(Y, \psi) $ which is unique up to
conjugacy. In particular, the canonical extension $(\ti{X},
\ti{\vpi})$ is such a homeomorphism extension.
\end{corollary}

If $(X, \vpi)$ is a dynamical system, then the map $ \al: C(X) \to
C(X),\ f \mapsto f \circ \vpi,$ is a $*$-endomorphism. $\al$ is a
$*$-automorphism iff $\vpi$ is a homeomorphism.  We can dualize the
preceding results as follows:

\begin{corollary} \label{c:al}
Given a dynamical system $(X, \vpi)$, there is a
minimal commutative C$^*$-algebra $C(\ti{X})$ with $*$ automorphism
$\ti{\al}$ admitting an embedding $\iota: C(X) \hookrightarrow
C(\ti{X})$ such that $\ti{\al}\circ \iota = \iota \circ \al .$
Furthermore, this commutative C$^*$-algebra is unique up to
isomorphism.
\end{corollary}

\begin{proof} Consider the inductive limit
\[ C(X) \stackrel{\al}{\to} C(X) \stackrel{\al}{\to} C(X)
\stackrel{\al}{\to} \dots .\] The inductive limit is a
C$^*$-algebra, $C(Y)$ containing $C(X)$ as a subalgebra, and $C(Y)$
admits an automorphism, $\be$ satisfying $\be(f) = \al(f)$ for $f
\in C(X)$.

But, with $(\ti{X}, \ti{\vpi})$ the minimal homeomorphism extension
of $(X, \vpi),$ and viewing $C(X) \hookrightarrow C(\ti{X})$, we can
consider the inductive limit
\[ C(X) \stackrel{id}{\to} \ti{\al}^{-1}(C(X)) \stackrel{id}{\to}
\ti{\al}^{-2}(C(X)) \dots .\]

The two inductive limits are isometrically isomorphic, as we have
the commutative diagram
\[
\begin{CD}
C(X) @>\ti{\al}^{-n}>>  \ti{\al}^{-n}(C(X)) \\
@V{\al}VV             @V{id}VV  \\
C(X) @>{\ti{\al}^{-(n+1)}}>> \ti{\al}^{-(n+1)}(C(X)))
\end{CD}
\]
Thus, we may identify $Y$ with $\ti{X}$, and we have the relation
\[ \ti{\al}^{-(n+1)}\al(f) =  \ti{\al}^{-n}(f), \ f \in C(X),\ n \in \bbZ^+, \]
hence
\[ \al(f) = \ti{\al}(f), \text{ or } \ti{\al}\circ \iota = \iota \circ \al
\] if we denote the embedding of $C(X) $ into $C(\ti{X})$ by
$\iota$.

\end{proof}

\begin{definition} \label{d:periodic}
Given a dynamical system $(X, \vpi)$, a point $x \in X$ is
\emph{periodic} if, for some $n \in \bbN,\ n \geq 1, \ \vpi^n(x) =
x$.  If $n$ is the smallest integer with this property, we say that
$x$ is periodic of period $n$.  If $x$ is not periodic, we say $x$
is \emph{aperiodic}. If for some $m \in \bbN,\ \vpi^m(x)$ is
periodic, then we say $x$ is \emph{eventually periodic}.
\end{definition}

\begin{remark} If $\vpi$ is a homeomorphism, then a point is
eventually periodic iff it is periodic; but if $\vpi$ is a
continuous surjection, it is possible to have a point $x$ which is
aperiodic and eventually periodic.
\end{remark}

\begin{lemma} \label{l:aperiodic} Let $(X, \vpi)$ be a dynamical system, and
$(\ti{X}, \ti{\vpi})$ its minimal homeomorphism extension.  Then a
point $\ti{x} = (x_1, x_2, \dots) \in \ti{X}$ is aperiodic iff for
any $n \in \bbN,\ x_n = x_m$ for at most finitely many $m \in \bbN
$, and $\ti{x}$ is periodic iff $x$ is periodic.
\end{lemma}

\begin{proof} $\ti{x}$ is periodic of period $p$ iff
$\ti{\vpi}^p(\ti{x}) = \ti{x},$ equivalently,
\begin{align*}
 (x_1, x_2, \dots) &= (\vpi^p(x_1), \vpi^p(x_2), \dots) \\
 &= (x_1, \dots, x_p, x_1, \dots, x_p, \dots)
 \end{align*}
which uses the relation that $\vpi^p(x_{p+j}) = x_j, \ j \in \bbN$.

 This shows that if $\ti{x}$ is periodic, the
coordinates of $\ti{x}$ form a periodic sequence; the converse is
also clear.
\end{proof}

\begin{definition} \label{d:dyndefs}
\begin{enumerate}
\item \label{i:toptrans} Recall a dynamical system $(X, \vpi)$ is \emph{topologically
transitive} if for any nonempty open set $\sO \subset X, \
\cup_{n=0}^{\infty} \vpi^{-n}\sO = X $.

\item \label{i:minimal} A dynamical system $(X, \vpi)$ is
\emph{minimal} if there is no proper, closed subset $Z \subset X$
such that $\vpi(Z) = Z $.

\item \label{i:recurr} A point $x$ in a dynamical system $(X, \vpi)$
is \emph{recurrent} if there is a subsequence $\{n_i\}$ of $\bbN$
such that $\vpi^{n_i}(x) \to x$.

\end{enumerate}
\end{definition}

\begin{remark} There should be no confusion between the two distinct
uses of \emph{minimal}.
\end{remark}

\begin{theorem} \label{t:extprop} Let $(X, \vpi)$ be a dynamical
system, and $(\ti{X}, \ti{\vpi})$ the minimal homeomorphism
extension.
\begin{enumerate}
\item $X$ is metrizable iff $\ti{X}$ is metrizable.

\item $(X, \vpi)$ is topologically transitive iff $(\ti{X},
\ti{\vpi})$ is topologically transitive.

\item $(X, \vpi)$ has a dense set of periodic points iff the same is
true of $(\ti{X}, \ti{\vpi})$.

\item $(X, \vpi)$ is a minimal dynamical system iff the minimal
homeomorphism extension has the same property.

\item The recurrent points
in $X$ are dense iff the recurrent points in $\ti{X}$ are dense.

\end{enumerate}

\end{theorem}

\begin{proof}
(1) is routine.

(2) Let $(\ti{X}, \ti{\vpi})$ be topologically transitive, and
$\emptyset \neq \sO \subset X.$  Then $\ti{\sO} := p^{-1}(\sO)$ is
nonempty in $\ti{X}$, so by assumption $\ti{X} = \cup_{n=0}^{\infty}
\ti{\vpi}^{-n}(\ti{\sO}) $.  Let $x \in X$ and $\ti{x} \in
p^{-1}(x)$. Then there exists $n$ such that $\ti{x} \in
\ti{\vpi}^{-n}(\ti{y}), \ \ti{y} \in \ti{\sO}$. So $\ti{y} =
\ti{\vpi}^n(\ti{x}),$ and so $p(\ti{y}) = p(\ti{\vpi}^n(\ti{x}) =
\vpi^n \circ p (\ti{x}) = \vpi^n(x)$.  Thus $x \in
\vpi^{-n}(p(\ti{y}) \subset \vpi^{-n}(\sO)$.

For the other direction, by Corollary~\ref{c:uniqminext} we can
assume, without loss of generality, that $(\ti{X}, \ti{\vpi})$ is
the canonical minimal homeomorphism extension of $(X, \vpi)$. The
basic open sets in $\ti{X}$ have the form
\[ \ti{\sO} = \ti{X} \cap [\sO_1 \times \dots \times \sO_N \times
\Pi_{n=N+1}^{\infty} X_n] \] for some $N \in \bbN, \sO_1, \dots
\sO_N $ open sets in $X$, and where $X_n = X$ for all $n>N$.

If $\ti{\sO}$ is nonempty, there is a point $x \in \sO_N$ such that
$\vpi^j(x) \in \sO_{N-j},\ j = 1, \dots, N-1.$  Hence by continuity
of $\vpi$ there is a neighborhood $U \subset \sO_N$ such that
$\vpi^j(U) \subset \sO_{N-j},\ j = 1, \dots, N-1$.

Let $\ti{x} \in \ti{X}$ be arbitrary, $\ti{x} = (x_1, \dots, x_N,
\dots)$.  By the topological transitivity of $(X, \vpi)$ we can find
$n \in \bbN$ such that $\vpi^n(x_N) \in U.$ Thus,
\begin{align*}
\ti{\vpi}^n(\ti{x}) &= (\vpi^n(x_1), \dots, \vpi^n(x_N), \dots) \\
    &\in \vpi^{N-1}(U)\times \dots \vpi(U) \times U \times
    \Pi_{n=N+1}^{\infty} X_n\\
    &\in \sO_1 \times \dots \times \sO_{N-1} \times \sO_N \times
    \Pi_{n=N+1}^{\infty}X_n
\end{align*}
so that $\ti{\vpi}^n(\ti{x}) \in \ti{\sO}$ which finishes the proof.

(3) If the periodic points are dense in $\ti{X}$, let $x \in X$ and
$\ti{x} \in p^{-1}(x)$. Then there is a net $\{\ti{y}_n\} \subset
\ti{X}$ converging to $x$. By Lemma~\ref{l:aperiodic}, $p(\ti{y}_n)$
is periodic in $X$, and converges to $x$.

For the converse, note that if $x \in X$ is periodic, say of period
$n$, then there is a point $\ti{x} \in \ti{X}, \ p(\ti{x}) = x$ with
$\ti{x}$ periodic of period $n$. Indeed, if $x$ has orbit $x,
\vpi(x), \dots, \vpi^{n-1}(x),$ then, setting $x_j =
\vpi^{n+1-j}(x),\ j = 1,\dots, n, $ take $\ti{x}$ to be the point
with coordinates
\[ \ti{x} = (x_1, x_2, \dots, x_n, x_1, x_2, \dots, x_n, \dots) .\]

If $\ti{\sO}$ is a basic open set in $\ti{X}$, we use the argument
in (2) to find an integer $N$ and an open set $U$ as in (2). Let $y
\in U$ be periodic, and set $x = \vpi^{N-1}(y)$. The point $\ti{x}
\in \ti{X}$ is periodic and belongs to $\ti{\sO}$.

(4) Assume $(\ti{X}, \ti{\vpi}) $ is a minimal dynamical system, and
$Y \subset X$ a nonempty closed, $\vpi$-invariant subset. Then
$p^{-1}(Y)$ is a nonempty closed invariant subset of $\ti{X}$, so
$p^{-1}(Y) = \ti{X}$.  Thus, $Y = X$, and so $(X, \vpi)$ is minimal.

Conversely, assume $(X, \vpi)$ be a minimal dynamical system, and
let $\sY = \{ Y_i \}_{i \in I}$ be a maximal chain of closed
invariant subsets of $\ti{X}$, ordered by inclusion. Then
\[ Y = \cap_{i \in I} Y_i \]
is the minimal element of the chain, hence $Y$ has no proper
invariant subset. As $Y \neq \emptyset,\ p(Y)$ is a nonempty
invariant subset of $X$, so $p(Y) = X$. Taking $\psi = \vpi|_Y,$ and
$ q = p|_Y$, we have that $(Y, \psi)$ is a homeomorphism extension
of $(X, \vpi)$. While we do not know \textit{a priori} that $(Y,
\psi)$ is a minimal homeomorphism extension of $(X, \vpi)$, if $(Y,
\psi)$ is not a minimal homeomorphism extension, there is an
intermediate extension $(Z, \si)$, as in the proof of
Lemma~\ref{l:minext}. As $\psi$ is a minimal homeomorphism on $Y$,
$\si$ is a minimal homeomorphism on $Z$.  It follows that the
minimal homeomorphism extension of $(X, \vpi)$ which lies between
$X$ and $Y$ is necessarily a minimal homeomorphism.

Since $(\ti{X}, \ti{\vpi})$ was the canonical minimal extension, and
since any two minimal extensions are conjugate, it follows that the
dynamical system $(\ti{X}, \ti{\vpi})$ is minimal.

(5) If the recurrent points in $\ti{X}$ are dense, let $U$ be any
nonempty open set in $X$. Then there exists $\ti{y} \in p^{-1}(U)$
which is recurrent.  But then $y := p(\ti{y}) \in U$ is recurrent.

Now assume $(X, \vpi)$ has a dense set of recurrent. First we show
that if $x \in X$ is recurrent, there is $\ti{x} \in p^{-1}(x)$
which is recurrent. So, let $x \in X$ be recurrent, and let $\ti{x}
= (x_1, x_2, \dots) \in \ti{X}$ be such that $p(\ti{x}) = x$ (so $x
= x_1$). By the compactness of $X$ and a standard diagonalization
argument, there is a subsequence $\{ n_j\}$ of $\bbN$ and $y_i \in
X,\ i = 1, 2, \dots$, such that
\[ \lim_j \vpi^{n_j}(x_i) = y_i,\ i = 1, 2, \dots \]
and $y_1 = x_1$. Since $\vpi(x_{i+1}) = x_i$, the same relation
holds for the $y_i$, and hence $\ti{y}:= (x_1, y_2, y_3, \dots) \in
\ti{X}$. Since
\[ \lim_j \vpi^{n_j}(y_i)= \lim_j \vpi^{n_j-i+1}(x_1) = \lim_j
\vpi^{n_j}(x_i) = y_i, \] $i = 1, 2, \dots$, this shows that $\ti{y}
\in p^{-1}(x)$ is recurrent.

Now, let $\ti{\sO} \subset \ti{X}$ be a basic open set, and let $U$
be an open set in $X$ and $N \in \bbN$ be as in the proof of (2).
Let $x_N \in U$ be recurrent; by the above assertion we can find
$\ti{x} = (x_1, \dots, x_N, \dots)$ which is recurrent in $X$, and
by construction $\ti{x}$ lies in $\ti{\sO}$.

\end{proof}

\section{Representations of Semicrossed Products} \label{s:semi}
For the moment we will take an abstract approach:  Let $(X, \vpi)$
be a dynamical system, and consider the algebra generated by $C(X)$
and a symbol $U$, where $U$ satisfies the relation
\[ (\ddag) \quad fU = U(f\circ \vpi),\ f \in C(X) .\]
The elements $F$ of this algebra can be viewed as noncommutative
polynomials in $U$,
\[ F = \sum_{n=0}^{N}  U^n\,f_n,\ f_n \in C(X), N \in \bbN. \]
Let us call this algebra $\sA_0 $.

In \cite{Pet.Sem} we formed the Banach Algebra $\ell_1(\sA_0)$ by
providing a norm to elements $F$ as above as $||F||_1 = \sum_{n=0}^N
||f_n||$ and then completing $\sA_0$ in this norm. Either approach
yields the same semicrossed product.

By a representation of $\sA_0$ we will mean a homomorphism of $\sA_0
$ into the bounded operators on a Hilbert space, which is a
$*$-representation when restricted to $C(X)$, viewed as a subalgebra
of $\sA_0 $, and such that $\pi(U)$ is an isometry.

Fix a point $x \in X$ and, for convenience, set $x_1 = x,\ x_2 =
\vpi(x),\ x_3 = \vpi^2(x), \ \dots$. Define a representation $\pi_x$
of $\sA_0$ on $\ell^2(\bbN)$ by
\[ \pi_x(f)(z_1, z_2, \dots) = (f(x_1)z_1, f(x_2)z_2, \dots), \]
with $(z_n)_{n=1}^{\infty} \in \ell^2(\bbN)$ and
\[ \pi_x(U)(z_1, z_2, \dots) = (0, z_1, z_2, \dots) .\]

Observe this is a representation of $\sA_0$ since

\[ \pi_x(fU)(z_1, z_2, \dots) = (0, f(x_2)z_1, f(x_3)z_2, \dots) \]
and
\begin{align*}
 \pi_x(Uf\circ \vpi)(z_1, z_2, \dots) &= \pi_x(U)(f\circ \vpi(x_1) z_1,
f\circ \vpi(x_2)z_2, \dots) \\
&= (0, f\circ \vpi(x_1) z_1, f\circ \vpi(x_2)z_2, \dots) \\
&= (0,  f(x_2)z_1, f(x_3)z_2, \dots)
\end{align*}

Let $(\ti{X}, \ti{\vpi})$ be the canonical homeomorphism extension.
(cf definition~\ref{d:homeoext} and Corollary~\ref{c:uniqminext}.)
We will consider $\sA_0$ as embedded in $\ti{\sA}_0,$ where
$\ti{\sA}_0$ is the algebra generated by $C(\ti{X})$ and $\ti{U}$,
satisfying the same relation ($\ddag$).  Let $\ti{x} \in \ti{X}$ and
set $x = p(\ti{x})$ where $p: \ti{X} \to X$ is the map in diagram
(\dag).

For $f \in C(X),$ let $\ti{f} \in C(\ti{X}),\ \ti{f} = f\circ p, $
and for $F = \sum_{n=0}^N U^n f_n,\ f_n \in C(X)$, let $\ti{F} =
\sum_{n=0}^N \ti{U}^n \ti{f}_n.$ Observe that
\[ \pi_{\ti{x}}(\ti{F}) = \pi_x(F) .\]

\subsection{Nest Representations}
For nonselfadjoint operator algebras, the representations which can
play the role of the primitive representations in the case of
C$^*$-algebras are the \emph{nest representations}.  Recall, a
representation $\pi$ of an algebra $\sA$ on a Hilbert space $\sH$ is
a nest representations if the lattice of subspaces invariant under
$\pi$ is linearly ordered.

Let $(X, \vpi)$ be a dynamical system with $\vpi$ a homeomorphism,
and $x$ a point  in  $X $ which is aperiodic. This means
$\ti{\vpi}^n(\ti{x}) \neq \ti{x}$ for all $ n \geq 1$.

\begin{lemma} \label{l:masa}  The weak closure of $\pi(C(X))$ is a
masa in $\sB(\sH).$
\end{lemma}

\begin{proof}  It is enough to show that the operator $e_n$ belongs
to the weak closure, where $e_n$ is the multiplication operator
which is $1$ in the $n^{th}$ coordinate and zero elsewhere.  We can
find $f_m \in C(X)$ satisfying
\[ f_m(x_j) =
\begin{cases}
1 &\text{ for } j = n \\
0 &\text{ for } j \neq n,\ j \leq m
\end{cases}
\]
and $f_m$ is real-valued, $ 0 \leq f_m \leq 1.$  Indeed, this
follows from the Tietze Extension Theorem.

As $\pi_x(f_m) \to e_n$ weakly, we have $e_n$ in the weak closure of
$\pi_x(C(X)$, and we are done.
\end{proof}

\begin{proposition} \label{p:nestrep} Let $(X, \vpi)$ be a dynamical system.
  If $x \in X $ is aperiodic, then $\pi_x $ is
a nest representation.
\end{proposition}

\begin{proof} Since the weak closure of $\pi(C(X))$ is a masa
(Lemma~\ref{l:masa}), the closed subspaces $\sS$ of $\ell^2(\bbN)$
invariant under $\pi(C(X))$ are the vectors $\ti{z} \in
\ell^2(\bbN)$ which are supported on a given subset of $\bbN$.  If
such a subspace is also invariant under $\pi(U)$ then it has the
form
\[ \sS = \{ \ti{z} \in \ell^2(\bbN): z_n = 0 \text{ for } n \leq N
\}
\] for some $N \in \bbN$. But then the subspaces $\sS$ are nested.
\end{proof}

To periodic points we can associate another class of nest
representations. Let $\ti{x} \in \ti{X}$ be periodic of period $N$,
so $ \ti{x} = (x_1, x_2, \dots)$ with $x_{i+N} = x_i$ for $i \in
\bbN$. Let $\pi: \sA \to \sB(\ell^2(N))$ by $\pi(f)(z_1, \dots, z_N)
= (f(x_1)z_1, \dots, f(x_N)z_N) $ and $\pi(U)(z_1, \dots, z_N) =
(z_N, z_1, \dots, z_{N-1})$.

For C$^*$ crossed products $\sB := C(X)\rtimes_{\psi} \bbZ $ where
$\psi$ is a homeomorphism, we have the representations $\Pi_x$ and
$\Pi_{y, \la} $ for $x$ aperiodic, $y$ periodic, and $\la \in \bbT$
given as follows: $\Pi_x$ acts on $\ell^2(\bbZ)$, where $\Pi_x(U)$
is the bilateral shift (to the right), and
\[ \Pi_x(f)(\xi_n) = (f(\psi^n(x))\xi_n) \quad n \in \bbZ, \ f \in
C(X). \]

$\Pi_{y, \la}$ acts on the finite dimensional space $\ell^2(p)$,
where $p$ is the period of the orbit of $y$. $\Pi_{y, \la}(U)$ is a
cyclic permutation along the (finite) orbit of $y$ composed with
multiplication by $\la$, and $\Pi_{y, \la}(f)$ acts like $\Pi_x(f)$
along the orbit of $y$.  These representations correspond to the
pure state extensions of the states on $C(X), \ f \to f(x)$ in the
cases where $x$ is aperiodic or periodic, respectively, and so are
irreducible.  However, not all irreducible representations of $\sB$
need be of this form. Nevertheless, Tomiyama has shown:

\begin{proposition} \label{p:tomi}
Every ideal of $\sB$ is the intersection of those ideals of the form
$ker(\Pi_x)$ and $ker(\Pi_{y, \la})$ ($x$ aperiodic, $y$ periodic,
$\la \ in \bbT$) which contain it.
\end{proposition}
This is Proposition 4.1 of \cite{Tom.Int}.

\begin{corollary} \label{c:tomi}
If $(Y, \si)$ is a dynamical system with $\si$ a homeomorphism, then
for $F \in C(Y)\rtimes_{\si}\bbZ$,
\[  ||F|| =  \max\{ A, B\}\]
where
\[ A = \sup \{ ||\Pi_x(F)||: x \text{ aperiodic}\} \]
and
\[ B = \sup\{ ||\Pi_{y, \la}(F)||: y \text{ periodic},\ \la \in \bbT
.\]
\end{corollary}

\begin{proof} Denote by $||\cdot||$ the crossed product norm, and by
$||\cdot||_{*}$ the norm defined in the statement of the corollary.
Let $I$ be the ideal in $C(Y)\rtimes_{\si}\bbZ$ of all $F$ with
$||F||_{*} = 0.$ Every ideal of the form $ker(\Pi_x)$ ($x$
aperiodic) and $ker(\Pi_{y, \la})$ ($y$ periodic,\ $\la \in \bbT$)
contains $I$.  Since the zero ideal also this property, it follows
from Proposition~\ref{p:tomi} that $I = (0)$.
\end{proof}

\begin{lemma} \label{l:norm1} Let $(X, \vpi)$ be a dynamical system with $\vpi$
a homeomorphism, and let $y \in X$ be periodic.  For $F \in
C(X)\rtimes_{\vpi} \bbZ^+,$ we have $||\pi_y(F)|| \geq \sup_{\la \in
\bbT} ||\Pi_{y, \la}(F)||$.
\end{lemma}

\begin{proof} Since any $F \in C(X)\rtimes_{\vpi} \bbZ^+$ can be
approximated by elements with finitely many nonzero Fourier
coefficients, we can assume $F$ has this property.  Let $y$ have
period $p$, and we can assume $F = \sum_{n=0}^{kp} U^nf_n, \ f_n \in
C(X)$, for some $k \in \bbZ^+$.

Let $\xi = (\xi_1, \dots, \xi_p) \in \bbC^p$ be any vector of norm
$1$, and fix $\la \in \bbT.$  For $N \in \bbN $ define a vector
$\eta \in \ell^2(\bbN)$ of norm $1$ by
\[ \eta = (\eta_1, \dots, \eta_{Np}, 0, 0, \dots) \text{ where } \eta_{i+jp} =
\la^{N-j} \xi_i/\sqrt{N},\ i = 1, \dots p,\ j = 0, \dots N-1 .\]

Now, for $k \leq j < N, $
\[ <\pi_y(F)\eta, e_{i+jp}> = \la^{j-k}/\sqrt{N}<\Pi_{y, \la}(F)\xi, e^p_i>
\] where $e_n$ resp. $e^p_n$ are standard basis vectors in
$\ell^2(\bbN)$, resp., in $\bbC^p$. Thus, if $N/k$ is large, it
follows that $||\pi_y(F)\eta||$ is close to $||\Pi_{y,
\la}(F)\xi||$.  This proves the lemma.
\end{proof}

\begin{lemma} \label{l:norm2} Let $(X, \vpi)$ be a dynamical system with $\vpi$
a homeomorphism, and $F \in C(X)\rtimes_{\vpi} \bbZ^+$. For any $x
\in X$,
\[ ||\Pi_x(F)|| = \sup \{ ||\pi_y(F)||:\ y \in \text{Orbit}(x) \} .\]
\end{lemma}

\begin{proof}
Given $\epsilon > 0,$ there is a vector $\xi \in \ell^2{\bbZ},\ \xi
= (\xi_n)_{n \in \bbZ}$ with only finitely many $\xi_n \neq 0$, and
such that
\[ ||\Pi_x(F)|| \leq ||\Pi_x(F)\xi|| + \epsilon .\]
Suppose $\xi_n = 0 $ for $n < -N,$ for some $N \in \bbZ^+ $. Let $y
= \vpi^{-N}(x),$ and define a vector $\eta \in \ell^2(\bbN)$ by:
$\eta_j = \xi_{j-N-1},\ j = 1, 2, \dots.$  Then $||\eta|| = 1,$ and
\[ ||\pi_y(F)\eta|| = ||\Pi_x(F)\xi|| .\]
The lemma now follows.
\end{proof}

\begin{corollary} Let $(X, \vpi)$ be a dynamical system, and
$(\ti{X} , \ti{\vpi})$ a minimal homeomorphism extension, with $p:
\ti{X} \to X$ the continuous surjection for which the diagram
($\dag$) commutes. Let $F \in C(X)\rtimes_{\vpi} \bbZ^+$, and let
$\ti{x} \in \ti{X}.$ Then
\[ ||\Pi_{\ti{x}}(\ti{F}|| = \sup\{ ||\pi_y(F)||:\ y = p(\ti{y}),
\text{ for } \ti{y} \in \text{Orbit}(\ti{x}) \} .\]
\end{corollary}

\begin{proof} Observe that for any $\ti{y} \in \ti{X}, $ and $\xi
\in \ell^2(\bbN),$
\[ \pi_{\ti{y}}(\ti{F})\xi = \pi_y(F)\xi .\]

Now apply Lemma~\ref{l:norm2}.
\end{proof}

\begin{definition} For a dynamical system $(X, \vpi)$ ($\vpi$ not
necessarily a homeomorphism), a periodic point $y \in X$ and $\la
\in \bbT$, we define $\pi_{y, \la}$ exactly like $\Pi_{y, \la}$ in
the case where $\vpi$ is a homeomorphism.
\end{definition}

\begin{remark} Since $\Pi_{y, \la}$ is irreducible, the same is true
for $\pi_{y, \la}$, and in particular $\pi_{y, \la}$ is a nest
representation.
\end{remark}

\begin{corollary} Let $(X, \vpi)$ be a dynamical system, $F \in
C(X)\rtimes_{\vpi} \bbZ^+$.  Then\[  ||F|| =  \max\{ A, B\}\] where
\[ A = \sup \{ ||\pi_x(F)||: x \text{ aperiodic}\} \]
and
\[ B = \sup\{ ||\pi_{y, \la}(F)||: y \text{ periodic},\ \la \in \bbT
.\]
\end{corollary}

\begin{proof} Note the constant "$A$" is the same as in Corollary~\ref{c:tomi} ,
and by Lemma~\ref{l:norm2} the constant "$B$" is the same as in
Corolary~\ref{c:tomi}.  For $y \in X$ periodic and $\la \in \bbT$,
we have by Lemma~\ref{l:norm1}
\[ \sup_{\la \in \bbT} ||\pi_{y, \la}(F)|| = \sup_{\la \in \bbT}
||Pi_{\ti{y}, \la}(\ti{F})|| \leq ||\pi_y(F)|| \leq ||\Pi_y(F)|| \]
where $\ti{y} \in \ti{X}$ is periodic and $ p(\ti{y}) = y$.

By Cor. II.8 of \cite{Pet.Sem}, $||F|| = \sup_{x \in X}
||\pi_x(F)||$.  Thus, denoting the $\sup\{A, B\}$ by $||F||_*$, by
Corollary~\ref{c:tomi} it follows that $||F||_* $ is the norm of $F$
in the crossed product $C(\ti{X})\rtimes_{\ti{\vpi}}\bbZ$.  On the
other hand,
\[ ||F||_* \leq \sup_{x\in X}||\pi_x(F)|| = ||F|| \leq \sup_{\ti{x}
\in \ti{X}} ||\Pi_{\ti{x}}(\ti{F})||\] and the last term is
dominated by the norm of $\ti{F}$ in the crossed product, since the
norm there is given by the supremum over all covariant
representations.
\end{proof}

\begin{theorem} \label{t:nestrep} Let $(X, \vpi)$ be a dynamical
system, $F \in C(X)\rtimes_{\vpi}\bbZ^+$.  Then
\[ ||F|| = \sup\{ ||\pi(F)||:\ \pi \text{is an isometric covariant nest
representation} \} .\]
\end{theorem}

\begin{proof} Indeed, we have found a subclass of the isometric
covariant nest representations, namely the $\pi_{y, \la}$ ($y$
periodic, $\la \in \bbT$) and $\pi_x$ ($x$ aperiodic) which yield
$||F||$.
\end{proof}

From the above, we obtain the following
\begin{theorem} \label{t:compisomemb}
Let $(X, \vpi)$ be a dynamical system, and $(\ti{X}, \ti{\vpi})$ its
minimal homeomorphism extension.  Then the embedding of the
semicrossed product $C(X) \rtimes_{\vpi} \bbZ^+ \hookrightarrow
C(\ti{X})\rtimes_{\ti{\vpi}} \bbZ $ into the crossed product is a
completely isometric isomorphism.
\end{theorem}

\begin{corollary} With notation as above, the semicrossed product
$C(X)\rtimes_{\vpi} \bbZ^+$ is semisimple iff
$C(\ti{X})\rtimes_{\ti{\vpi}} \bbZ^+ $ is semisimple.
\end{corollary}

\begin{proof} This follows from part (5) of Theorem~\ref{t:extprop}
and the main result of \cite{DKM.Jac}.
\end{proof}

If the crossed product is a simple C$^*$-algebra, the crossed
product is necessarily the C$^*$-envelope.  However, as we will now
show, it is always the case that the crossed product is the
C$^*$-envelope, even if it is not simple.

\begin{lemma} \label{l:endo}  Let $(X, \vpi)$ be a dynamical system,
and $C(X) \rtimes_{\vpi} \bbZ^+$ the associated semicrossed product.
Then the endomorphism $\al$ of $C(X),\ \al(f) = f \circ \vpi,$
extends to an endomorphism, again denoted by $\al$, of the
semicrossed product.
\end{lemma}

\begin{proof} Embed $C(X) \rtimes_{\vpi} \bbZ^+ \hookrightarrow
C(\ti{X})\rtimes_{\ti{\vpi}} \bbZ $. The element $U$, which is an
isometry in the semicrossed product, embeds to a unitary in the
crossed product, and one can define an automorphism $\ti{\al}$ on
the crossed product, which extends the automorphism, also denoted by
$\ti{\al}$ of $C(\ti{X}),\ \ti{\al}(f) = f\circ \ti{\vpi}.$

This is as follows: in $C(\ti{X})\rtimes_{\ti{\vpi}} \bbZ $ one has
$U^*fU = f\circ \ti{\vpi}$. For $F$ an element of the crossed
product, define $\ti{\al}(F) = U^*FU$.  Note that, if $\{f_n\}$ are
the Fourier coefficients of $F$, then $\{f_n\circ \ti{\vpi}\}$ are
the Fourier coefficients of $\ti{\al}(F)$.

In particular, if $F$ belongs to the semicrossed product, and so its
Fourier coefficients belong to the subalgebra $C(X) \hookrightarrow
C(\ti{X}),$ then the Fourier coefficients of $\ti{\al}(F)$ also
belong to the subalgebra $C(X)$, since for $f \in C(X),\ \ti{\al}(f)
= \al(f)$.  Thus, if we denote this map of $C(X) \rtimes_{\vpi}
\bbZ^+$ by $\al$, it is an endomorphism of the semicrossed product
extending the endomorphism $\al$ of $C(X)$.
\end{proof}

\begin{lemma} \label{l:densesub} Let $(X, \vpi)$ be a dynamical
system, and embed $ \sA := $\mbox{$C(X) \rtimes_{\vpi} \bbZ^+ $}$
\hookrightarrow C(\ti{X})\rtimes_{\ti{\vpi}} \bbZ $.  Then
\[ \cup_{k=0}^{\infty} \ti{\al}^{-k}(\sA) \subset
C(\ti{X})\rtimes_{\ti{\vpi}} \bbZ  \text{ is a dense subalgebra}.\]
\end{lemma}

\begin{proof} It follows from the proof of Corollary~\ref{c:al}
that, viewing $C(X)$ as embedded in $C(\ti{X})$, that
$\cup_{n=0}^{\infty} \ti{\al}^{-n}(C(X))$ is a dense subalgebra of
$C(\ti{X})$.

Given $F \in C(\ti{X})\rtimes_{\ti{\vpi}} \bbZ$ and $\epsilon > 0$,
there is $G \in C(\ti{X})\rtimes_{\ti{\vpi}} \bbZ$ with finitly many
 nonzero Fourier coefficients, say $G = \sum_{n=0}^N U^n g_n,$
 with $||F - G|| < \epsilon$. By the first paragraph, for each $g_n$
 there is an $h_n$ in the dense subalgebra of $C(\ti{X})$ with $||g_n - h_n|| <
 \frac{\epsilon}{N+1}$.  But if $H = \sum_{n=0}^{N}U^n h_n,$ we have $||F -
 H|| < 2\epsilon$, and $H \in \cup_{k=0}^{\infty}
 \ti{\al}^{-k}(\sA)$.

 \end{proof}

 \begin{theorem} \label{t:maint} Let $(X, \vpi)$ be a dynamical
 system, and $(\ti{X}, \ti{\vpi})$ its minimal homeomorphism
 extension. Then the C$^*$-envelope of the semicrossed product
 $C(X) \rtimes_{\vpi} \bbZ^+ $ is the crossed product
 $C(\ti{X})\rtimes_{\ti{\vpi}} \bbZ$.
 \end{theorem}

 \begin{proof} By Theorem~\ref{t:compisomemb} the embedding
 \[C(X) \rtimes_{\vpi} \bbZ^+ \hookrightarrow C(\ti{X})\rtimes_{\ti{\vpi}} \bbZ
 \] is completely isometric. Suppose there is a C$^*$-algebra $\sB$,
 a completely isometric embedding $\iota: C(X) \rtimes_{\vpi} \bbZ^+
 \to \sB$, and a surjective C$^*$-homomorphism
 $ q: C(\ti{X})\rtimes_{\ti{\vpi}} \bbZ \to \sB$.

 If $q$ is not an isomorphism, let $0 \neq F \in \text{ker}(q).$ Assume
 $||F|| = 1$.  By Lemma~\ref{l:densesub} there is an element $G =
 \sum_{n=0}^N U^n g_n$ with $g_n \in \cup_{n=0}^{\infty}
 \ti{\al}^{-n}(C(X)),$ viewing $C(X)$ as a subalgebra of
 $C(\ti{X})$, and such that $||F - G|| < \frac{1}{2}$.  In
 particular, there is $m \in \bbZ^+ $ such that $g_n \circ
 \ti{\vpi}^m \in C(X),\ 0 \leq n \leq N$.

 Now  $GU^m = \sum_{n=0}^N g_n \circ \ti{\vpi}^m \in C(X) \rtimes_{\vpi}
 \bbZ^+$ and $||GU^m|| = ||G|| > \frac{1}{2}$. On the other hand,
 since $q(FU^m) = q(F) q(U^m) = 0$ we have
 \[ ||q(GU^m)|| = ||q(GU^m - FU^m)|| \leq ||GU^m - FU^m|| \leq ||G -
 F|| < \frac{1}{2} .\]

 This contradiction shows that ker$(q) = (0)$, and hence that
 $C(\ti{X})\rtimes_{\ti{\vpi}} \bbZ$ is the C$^*$ envelope of the
 semicrossed product.

 \end{proof}

Finally, we make use of the relation between properties of dynamical
systems and their extensions to obtain

\begin{proposition} \label{p:semis}
Let $(X, \vpi)$ be a dynamical system. If the
C$^*$-envelope of the semicrossed product is a simple C$^*$-algebra,
then $C(X)\rtimes_{\vpi} \bbZ^+$ is semi-simple.
\end{proposition}

\begin{remark} The converse is false.
\end{remark}

\begin{proof} By Theorem~\ref{t:maint}, the C$^*$-envelope is
a crossed product, $C(\ti{X})\rtimes_{\ti{\vpi}}\bbZ$, where
$(\ti{X}, \ti{\vpi})$ is the (unique) minimal homeomorphism
extension of $(X, \vpi)$. As is well known (e.g. \cite{Ped.Cstar}
Proposition 7.9.6), the crossed product is simple if and only if the
dynamical system $(\ti{X}, \ti{\vpi})$ is minimal; i.e., every point
has a dense orbit. By Theorem~\ref{t:extprop}, this is equivalent to
the condition that $(X, \vpi)$ is minimal.  In particular, the
system $(X, \vpi)$ is recurrent; so by \cite{DKM.Jac} it follows
that the semicrossed product is semi-simple.
\end{proof}

\end{document}